\documentclass[a4paper,11pt]{amsart}
\usepackage{amssymb,amscd}
\usepackage{epsfig}
\usepackage[cp1251]{inputenc}
\usepackage[T2A]{fontenc}
\newtheorem{theorem}{Theorem}[section]
\newtheorem{proposition}[theorem]{Proposition}
\newtheorem{corollary}[theorem]{Corollary}
\newtheorem{lemma}[theorem]{Lemma}

\theoremstyle{definition}
\newtheorem{definition}[theorem]{Definition}
\newtheorem{example}[theorem]{Example}
\newtheorem{construction}[theorem]{Construction}
\theoremstyle{remark}
\newtheorem*{remark}{Remark}
\newtheorem*{acknowledgements}{Acknowledgements}

\numberwithin{equation}{section}

\def\C{\mathbb C}
\def\Q{\mathbb Q}
\def\R{\mathbb R}
\def\Z{\mathbb Z}

\def\sK{\mathcal K}
\def\sS{\mathcal S}

\def\phi{\varphi}

\def\ge{\geqslant}

\def\le{\leqslant}

\newcommand{\mdeg}{\mathop{\rm mdeg}}

\newcommand{\Tor}{\mathop{\rm Tor}\nolimits}
\def\rk{\mathop{\mathrm{rk}}}
\def\mrk{\mathop{\mathrm{mrk}}}

\def\hatzero{\hat 0}

\newcommand{\cs}{\mathbin{\#}}

\def\hrk{\mathop{\mathrm{hrk}}}
\def\trk{\mathop{\mathrm{trk}}}

\newcommand{\zk}{\mathcal Z_{\mathcal K}}
\def\zs{\mathcal Z_{\mathcal S}}

\def\zq{\mathcal Z_Q}

\def\cat{{\text{\sc cat}}}
\def\cgr{\text{\sc cgr}}
\def\cgg{\text{\sc cgg}}
\def\top{\text{\sc top}}
\def\colim{\mathop\mathrm{colim}\nolimits}
\def\lim{\mathop\mathrm{lim}\nolimits}

\newcommand{\mb}[1]{{\textbf {\textit#1}}}

\begin{document}

\title[Moment-angle complexes from simplicial posets]{Moment-angle complexes from\\simplicial posets}

\author{Zhi L\"u}
\author{Taras Panov}\thanks{The first author was
supported by the NSFC, grants~10671034 and~10931005. The second
author was supported by the Russian Foundation for Basic Research,
grants~10-01-92102-ЯФ and~09-01-00239; grants НШ-5413.2010.1 and
МД-2253.2011.1 from the President of Russia; and a grant from
Dmitri Zimin's `Dynasty' foundation. Both authors were also
supported by the SGST, grant~09DZ2272900.}
\address{Institute of Mathematics, School of Mathematical Sciences, Fudan University,
Shanghai, 200433, People's Republic of China}
\address{Department of Geometry and Topology, Faculty of Mathematics and Mechanics, Moscow State University,
Leninskie Gory, Moscow 119991, Russia\newline \emph{and}
\newline Institute for Theoretical and Experimental Physics, Moscow
117259, Russia}
\email{zlu@fudan.edu.cn}
\email{tpanov@mech.math.msu.su}


\begin{abstract}
We extend the construction of moment-angle complexes to simplicial
posets by associating a certain $T^m$-space $\zs$ to an arbitrary
simplicial poset $\sS$ on $m$ vertices. Face rings $\Z[\sS]$ of
simplicial posets generalise those of simplicial complexes, and give rise to
new classes of Gorenstein and Cohen--Macaulay rings. Our primary
motivation is to study the face rings $\Z[\sS]$ by topological
methods. The space $\zs$ has many important topological properties
of the original moment-angle complex $\zk$ associated to a
simplicial complex~$\sK$. In particular, we prove that the
integral cohomology algebra of $\zs$ is isomorphic to the
$\Tor$-algebra of the face ring $\Z[\sS]$. This leads directly to
a generalisation of Hochster's theorem, expressing the algebraic
Betti numbers of the ring $\Z[\sS]$ in terms of the homology of
full subposets in~$\sS$. Finally, we estimate the total amount of
homology of $\zs$ from below by proving the toral rank conjecture
for the moment-angle complexes~$\zs$.
\end{abstract}

\maketitle

\section{Introduction}
Simplicial posets describe the combinatorial structures underlying
``generalised simplicial complexes'' whose faces are still
simplices, but two faces are allowed to intersect in any
subcomplex of their boundary, rather than just in a single face.
These are also known as ``ideal triangulations'' in
low-dimensional topology, or as ``simplicial cell complexes''.

Simplicial posets also attract a lot of interest in algebraic
combinatorics and \emph{combinatorial commutative algebra}.
Stanley~\cite{stan91} introduced the \emph{face ring} $\Z[\sS]$ of
a simplicial poset $\sS$ as a quotient of certain graded
polynomial ring by a homogeneous ideal determined by the poset
relation in~$\sS$ (see Definition~\ref{deffr} below). The ring
$\Z[\sS]$ generalises the \emph{Stanley--Reisner face ring}
$\Z[\sK]$ of a simplicial complex~$\sK$. The rings~$\Z[\sS]$ have
remarkable algebraic and homological properties, albeit they are
much more complicated than the Stanley--Reisner rings~$\Z[\sK]$.
Unlike~$\Z[\sK]$, the ring~$\Z[\sS]$ is not generated in the
lowest positive degree. (For topological reasons it is convenient
to double the grading making it even; so that $\Z[\sK]$ is
generated by degree-two elements, but $\Z[\sS]$ is not.) Face
rings of simplicial posets were further studied by
Duval~\cite{duva97} and Maeda--Masuda--Panov \cite{ma-pa06},
\cite{m-m-p07}, among others. \emph{Cohen--Macaulay} and
\emph{Gorenstein*} face rings are particularly important; both
properties are topological, that is, depend only on the
topological type of the geometric realisation~$|\sS|$. Gorenstein*
simplicial posets also feature in toric topology, as combinatorial
structures associated to orbit quotients of \emph{torus manifolds}
with special cohomological properties~\cite{ma-pa06}.

Here we suggest an approach to studying the face rings of
simplicial posets by topological methods. We associate to $\sS$ a
certain space $\zs$, called the \emph{moment-angle complex}, which
is glued from products of discs and circles
(Definition~\ref{defmac}). The original moment-angle complex was
introduced by Buchstaber and Panov in~\cite{bu-pa99} as a
disc-circle decomposition of the Davis--Januszkiewicz universal
space $\zk$ associated to a simplicial
complex~$\sK$~\cite{da-ja91}; this decomposition was used in the
calculation of the cohomology ring of $\zk$ in terms of the face
ring of $\sK$~\cite{bu-pa99}.

We therefore continue here the unifying theme of toric topology
which links several aspects of equivariant topology to
combinatorial commutative algebra. Motivated by the categorical
constructions in toric topology~\cite{pa-ra08} we describe the
face ring $\Z[\sS]$ as the (inverse) limit of a certain diagram of
polynomial rings over the opposite face category $\cat^{op}(\sS)$
of~$\sS$ (Lemma~\ref{zslim}). This generalises the limit
description~\cite[(4.7)]{p-r-v04} for the Stanley--Reisner face
ring $\Z[\sK]$ of a simplicial complex, and leads to an important
functorial property (Proposition~\ref{funct}).

The face ring $\Z[\sS]$ of a simplicial poset $\sS$ with $m$
vertices is naturally an algebra over the polynomial ring
$\Z[v_1,\ldots,v_m]$. We show that the corresponding
\emph{$\Tor$-algebra} $\Tor_{\Z[v_1,\ldots,v_m]}(\Z[\sS],\Z)$ is
isomorphic to the integral cohomology ring of $\zs$ with the
appropriately refined grading (Theorem~\ref{hzs1}), therefore
extending the result of~\cite{bu-pa02} and~\cite{pano08l} to
simplicial posets. The Koszul complex calculating the $\Tor$
splits into the sum of subcomplexes corresponding to the
\emph{full subposets} $\mathcal S_\mb a$ of $\sS$; the cohomology
of these subcomplexes can be identified with the cellular
cohomology of $|\mathcal S_\mb a|$ after a shift of dimension.
This leads to a generalisation of Hochster's theorem calculating
the algebraic Betti numbers of~$\Z[\sS]$ (Corollary~\ref{gh}).

Recently a lot of work has been done on generalising the
construction of moment-angle complex $\zk=\zk(D^2,S^1)$ to pairs
of spaces $(X,W)$ different from $(D^2,S^1)$, and studying the
resulting spaces from the homotopy-theoretical perspective.
See~\cite{gr-th06} and~\cite{b-b-c-g07} for important advances in
this direction. Among examples of these ``generalised moment-angle
complexes'' we mention those corresponding to the pairs
$(D^1,S^0)$ (the \emph{real moment-angle complex} $\R\zk$),
$(\C,\C\setminus0)$ and $(\R,\R\setminus0)$ (the complex and real
\emph{coordinate subspace arrangement complements} respectively),
and $(\C P^\infty,pt)$ (the so-called \emph{Davis--Januszkiewicz
space}, whose cohomology is the face ring~$\Z[\sK]$),
see~\cite[Ch.~6]{bu-pa02}.

Here we follow a different route: instead of replacing the pair
$(D^2,S^1)$ in $\zk=\zk(D^2,S^1)$ by a different pair, we extend
the ``indexing structure'' from a simplicial complex $\sK$ to a
simplicial poset~$\sS$. One of the main reasons to keep the pair
$(D^2,S^1)$ intact is that the space $\zs=\zs(D^2,S^1)$ supports a
$T^m$-action, like the original moment-angle complex~$\zk$.
Moreover, if the dimension of $|\sS|$ is $n-1$, then there is
always an $(m-n)$-dimensional subtorus in $T^m$ acting on $\zs$
almost freely (Corollary~\ref{trkzs}). A choice of such subtorus
is equivalent to a choice of a linear system of parameters in the
$\Q$-face ring $\Q[\sS]$ (Theorem~\ref{zkquoq}). It has been shown
recently by Cao--L\"u~\cite{ca-lu09} and Ustinovsky~\cite{usti09}
that the total dimension of the rational cohomology of~$\zk$ is at
least $2^{m-n}$. Here we extend this result to $\zs$, thereby
settling Halperin's \emph{toral rank conjecture} for moment-angle
complexes corresponding to simplicial posets
(Corollary~\ref{trczs}).

%

We work over $\Z$ throughout most of the paper, as this is the
most natural coefficient ring from topologist's point of view. All
our statements are readily generalised to an arbitrary commutative
associative ring with unit.

There is a clash of terminology between combinatorialists and
homotopy theorists about using the term ``simplicial''. We do not
use simplicial homotopy theory in this paper, so that our
simplicial posets and simplicial cell complexes \emph{do not mean}
simplicial objects in the appropriate categories.

\begin{acknowledgements}
The authors are grateful to Peter Landweber for his careful
reading of the preliminary version of the manuscript and providing
a number of valuable comments.
\end{acknowledgements}

\section{Simplicial posets and their face rings}
A poset (partially ordered set) $\mathcal S$ with the order
relation $\le$ is called \emph{simplicial} if it has an initial
element $\hatzero$ and for each $\sigma\in\mathcal S$ the lower
segment $[\hatzero,\sigma]$ is a boolean lattice (the face poset
of a simplex). We assume all our posets to be finite and sometimes
refer to elements $\sigma\in\sS$ as \emph{simplices}. The
\emph{rank function} $|\cdot|$ on $\sS$ is defined by setting
$|\sigma|=k$ for $\sigma\in\mathcal S$ if $[\hatzero,\sigma]$ is
the face poset of a $(k-1)$-dimensional simplex. The rank of $\sS$
is the maximum of ranks of its elements, and the \emph{dimension}
of $\sS$ is its rank minus one. The \emph{vertices} of $\sS$ are
elements of rank one. We assume that $\sS$ has $m$ vertices and
denote the vertex set by $V(\sS)=[m]=\{1,\ldots,m\}$. Similarly we
denote by $V(\sigma)$ the vertex set of $\sigma$, that is the set
of $i$ with $i\le\sigma$.

The face poset of a simplicial complex is a simplicial poset, but
there are many simplicial posets that do not arise in this way. We
identify a simplicial complex with its face poset, thereby
regarding simplicial complexes as particular cases of simplicial
posets.

To each $\sigma\in\mathcal S$ we assign a geometric simplex
$\Delta^\sigma$ whose face poset is $[\hatzero,\sigma]$, and glue
these geometric simplices together according to the order relation
in $\mathcal S$. We get a regular cell complex in which the
closure of each cell is identified with a simplex preserving the
face structure, and all attaching maps are inclusions. We call it
a \emph{simplicial cell complex} and denote its underlying space
by $|\mathcal S|$.

Using a more formal categorical language, we consider the
\emph{face category} $\cat(\sS)$ whose objects are elements
$\sigma\in\sS$ and there is a morphism from $\sigma$ to $\tau$
whenever $\sigma\le\tau$. Then we may write
\[
  |\sS|=\colim \Delta^\sS,
\]
where $\Delta^\sS$ is a diagram (covariant functor) from
$\cat(\sS)$ sending every morphism $\sigma\le\tau$ to the
inclusion of geometric simplices
$\Delta^\sigma\hookrightarrow\Delta^\tau$, and the colimit is
taken in the category of (good) topological spaces.


For every simplicial poset $\sS$ there is the associated
simplicial complex $\sK_\sS$ on the same vertex set~$V(\sS)$,
whose simplices are the sets $V(\sigma)$, \ $\sigma\in\sS$. There
is a \emph{folding map} of simplicial posets
\begin{equation}\label{fold}
  \sS\longrightarrow\sK_\sS,\quad\sigma\mapsto V(\sigma).
\end{equation}
The corresponding geometric folding $|\sS|\to|\sK_\sS|$ is a
``branched combinatorial covering'' in the sense
of~\cite{bu-pa04-1}; it is the identity on the vertices, and every
simplex in $\sK_\sS$ is covered by a certain positive number of
simplices of~$\sS$.

For any two simplices $\sigma,\tau\in\mathcal S$, denote by
$\sigma\vee\tau$ the set of their least common upper bounds
(\emph{joins}), and by $\sigma\wedge\tau$ the set of their
greatest common lower bounds (\emph{meets}). Since $\mathcal S$ is
a simplicial poset, $\sigma\wedge\tau$ consists of a single
simplex whenever $\sigma\vee\tau$ is non-empty. It is easy to
observe that $\sS$ is a simplicial complex if and only if for any
$\sigma,\tau\in\sS$ the set $\sigma\vee\tau$ is either empty or
consists of a single simplex~\cite[Prop.~5.1]{m-m-p07}. In this
case $\sS$ coincides with $\sK_\sS$.

Now consider the graded polynomial ring $\Z[v_\sigma\colon
\sigma\in\mathcal S]$ with one generator $v_\sigma$ of degree
$\deg v_\sigma=2|\sigma|$ for every $\sigma\in\mathcal S$.

\begin{definition}[\cite{stan91}]\label{deffr}
The \emph{face ring} of a simplicial poset $\mathcal S$ is the
quotient
\[
  \Z[\mathcal S]:=
  \Z[v_\sigma\colon \sigma\in\mathcal S]\;/\;
  \mathcal I_{\mathcal S},
\]
where $\mathcal I_{\mathcal S}$ is the ideal generated by the
elements $v_{\hatzero}-1$ and
\begin{equation}\label{estre}
  v_\sigma v_\tau-v_{\sigma\wedge\tau}\cdot\!\!\!\sum_{\eta\in\sigma\vee\tau}\!\!\!v_\eta.
\end{equation}
The sum over the empty set is assumed to be zero, so we have
$v_\sigma v_\tau=0$ in $\Z[\sS]$ if $\sigma\vee\tau=\varnothing$.

The grading may be refined to a $\Z^m$-grading by setting $\mdeg
v_\sigma=2V(\sigma)$. Here $V(\sigma)$ is a subset of $[m]$, and
we identify such subsets $\mb a\subset[m]$ with vectors
in~$\{0,1\}^m\subset\Z^m$ in the standard way: the unit
coordinates of a vector correspond to the elements in a subset. In
particular, $\mdeg v_i=2\mb e_i$ (two times the $i$th basis
vector).
\end{definition}

\begin{remark}
The definition above extends the notion of the face ring of a
simplicial complex (also known as the \emph{Stanley--Reisner
ring}) to simplicial posets. In the case when $\mathcal S$ is a
simplicial complex we may rewrite~\eqref{estre} as $v_\sigma
v_\tau-v_{\sigma\wedge\tau}v_{\sigma\vee\tau}$ (because
$\sigma\vee\tau$ is either empty or consists of a single simplex),
and use the latter relation to express any $v_\sigma$ as
\[
  v_\sigma=\prod_{i\in V(\sigma)}v_i.
\]
The relations between the $v_i$ coming from~\eqref{estre} can now
be written as
\begin{equation}\label{screl}
v_{i_1}\cdots v_{i_k}=0\quad\text{if }\{i_1,\ldots,i_k\}\text{
does not span a simplex of }\sS.
\end{equation}
The face ring $\Z[\sS]$ is therefore isomorphic to the quotient of
the polynomial ring $\Z[v_1,\ldots,v_m]$ by~\eqref{screl}, where
$\deg v_i=2$. This is the standard way of describing the face ring
of a simplicial complex.
\end{remark}

\begin{example}\label{exsp}
\begin{figure}
\begin{picture}(0,0)%
\includegraphics{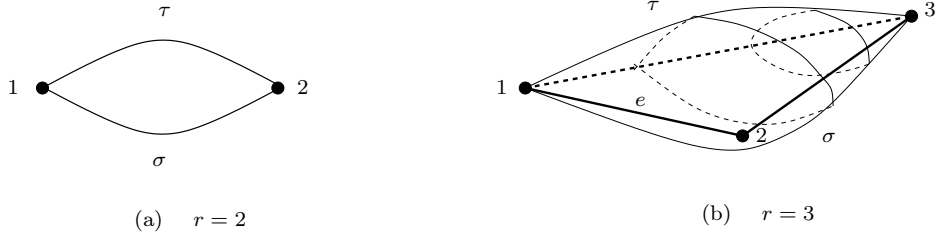}%
\end{picture}%
\setlength{\unitlength}{1973sp}%
\begingroup\makeatletter\ifx\SetFigFont\undefined%
\gdef\SetFigFont#1#2#3#4#5{%
  \reset@font\fontsize{#1}{#2pt}%
  \fontfamily{#3}\fontseries{#4}\fontshape{#5}%
  \selectfont}%
\fi\endgroup%
\begin{picture}(11769,2907)(676,-5149)
\put(10801,-4036){\makebox(0,0)[lb]{\smash{{\SetFigFont{8}{9.6}{\rmdefault}{\mddefault}{\updefault}
$\sigma$}}}}
\put(4276,-3436){\makebox(0,0)[lb]{\smash{{\SetFigFont{8}{9.6}{\rmdefault}{\mddefault}{\updefault}
$2$}}}}
\put(676,-3436){\makebox(0,0)[lb]{\smash{{\SetFigFont{8}{9.6}{\rmdefault}{\mddefault}{\updefault}
$1$}}}}
\put(2551,-2461){\makebox(0,0)[lb]{\smash{{\SetFigFont{8}{9.6}{\rmdefault}{\mddefault}{\updefault}
$\tau$}}}}
\put(2476,-4336){\makebox(0,0)[lb]{\smash{{\SetFigFont{8}{9.6}{\rmdefault}{\mddefault}{\updefault}
$\sigma$}}}}
\put(2251,-5086){\makebox(0,0)[lb]{\smash{{\SetFigFont{8}{9.6}{\rmdefault}{\mddefault}{\updefault}
(a)   \ \ \ $ r=2$}}}}
\put(9301,-5011){\makebox(0,0)[lb]{\smash{{\SetFigFont{8}{9.6}{\rmdefault}{\mddefault}{\updefault}
(b) \ \ \ $ r=3$}}}}
\put(6751,-3436){\makebox(0,0)[lb]{\smash{{\SetFigFont{8}{9.6}{\rmdefault}{\mddefault}{\updefault}
$1$}}}}
\put(12076,-2461){\makebox(0,0)[lb]{\smash{{\SetFigFont{8}{9.6}{\rmdefault}{\mddefault}{\updefault}
$3$}}}}
\put(9976,-4036){\makebox(0,0)[lb]{\smash{{\SetFigFont{8}{9.6}{\rmdefault}{\mddefault}{\updefault}
$2$}}}}
\put(8626,-2386){\makebox(0,0)[lb]{\smash{{\SetFigFont{8}{9.6}{\rmdefault}{\mddefault}{\updefault}
$\tau$}}}}
\put(8476,-3586){\makebox(0,0)[lb]{\smash{{\SetFigFont{8}{9.6}{\rmdefault}{\mddefault}{\updefault}
$e$}}}}
\end{picture}%
\centering
   \caption{Simplicial cell complexes.} \label{sccfig}
        \end{figure}

1. The simplicial cell complex shown in Fig.~\ref{sccfig}~(a) is
obtained by gluing two segments along their boundaries and has
rank~$2$. The vertices are $1,2$ and we denote the 1-dimensional
simplices by $\sigma$ and $\tau$. Then the face ring $\Z[\sS]$ is
the quotient of the graded polynomial ring
\[
  \Z[v_1,v_2,v_\sigma,v_\tau],\quad\deg v_1= \deg v_2=2,\quad
  \deg v_\sigma=\deg v_\tau=4
\]
by the two relations
\[
  v_1v_2=v_\sigma+v_\tau,\quad v_\sigma v_\tau=0.
\]

2. The simplicial cell complex in Fig.~\ref{sccfig}~(b) is
obtained by gluing two triangles along their boundaries and has
rank~$3$. The vertices are $1,2,3$ and we denote the 1-dimensional
simplices (edges) by $e$, $f$ and $g$, and the 2-dimensional
simplices by $\sigma$ and~$\tau$. The face ring $\Z[\sS]$ is
isomorphic to the quotient of the graded polynomial ring
\[
  \Z[v_1,v_2,v_3,v_\sigma,v_\tau],\quad\deg v_1= \deg v_2=\deg v_3=2,
  \quad \deg v_\sigma=\deg v_\tau=6
\]
by the two relations
\[
  v_1v_2v_3=v_\sigma+v_\tau,\quad v_\sigma v_\tau=0.
\]
The generators corresponding to the edges can be excluded because
of the relations $v_e=v_1v_2$, etc.
\end{example}

The following lemma gives another perspective on the algebraic
structure of the ring $\Z[\sS]$.

\begin{lemma}[{\cite[Lemma~5.4]{ma-pa06}}]\label{chain}
Every element of $\Z[\sS]$ can be uniquely written as a linear
combination of monomials
$v_{\tau_1}^{\alpha_1}v_{\tau_2}^{\alpha_2}\cdots
v_{\tau_k}^{\alpha_k}$ corresponding to chains of totally ordered
elements $\tau_1<\tau_2<\cdots<\tau_k$
of~$\sS\!\setminus\!\hatzero$.
\end{lemma}

In other words, the monomials
$v_{\tau_1}^{\alpha_1}v_{\tau_2}^{\alpha_2}\cdots
v_{\tau_k}^{\alpha_k}$ with $\tau_1<\tau_2<\cdots<\tau_k$
constitute a basis for the graded free abelian group~$\Z[\sS]$. We
refer to the expansion of an element $x\in\Z[\sS]$ in terms of
this basis as the \emph{chain decomposition} of~$x$. The proof of
the above lemma uses the \emph{straightening
relation}~\eqref{estre} inductively, which allows one to express a
product of two elements via products of elements in order. This
can be restated by saying that $\Z[\sS]$ is an example of an
\emph{algebra with straightening law} (see discussion
in~\cite[p.~323]{stan91}).

As was observed in~\cite[(4.7)]{p-r-v04}, the face ring $\Z[\sK]$
of a simplicial complex can be realised as the limit of a diagram
of polynomial algebras over $\cat^{op}(\sK)$. A similar
description exists for the face ring $\Z[\sS]$.

\begin{construction}[{$\Z[\sS]$ as limit}]
We consider the diagram (covariant functor) $\Z[\cdot]_\mathcal S$
from the opposite face category $\cat^{op}(\sS)$ to the category
$\cgr$ of commutative associative graded rings with unit. Its
value on $\sigma\in\sS$ is the polynomial ring
$\Z[\sigma]=\Z[v_i\colon i\in V(\sigma)]$, and its value on the
morphism $\sigma\le\tau$ is the surjection $\Z[\tau]\to\Z[\sigma]$
sending each $v_i$ with $i\notin V(\sigma)$ to zero.

\begin{lemma}\label{zslim}
We have
\[
  \Z[\sS]=\lim\Z[\cdot]_\mathcal S
\]
where the (inverse) limit is taken in the category $\cgr$.
\end{lemma}
\begin{proof}
We enumerate the elements of~$\sS$ so that the rank function does
not decrease, and proceed by induction. We therefore may assume
the statement is proved for a simplicial poset~$\mathcal T$, and
need to prove it for $\sS$ which is obtained from $\mathcal T$ by
adding one element~$\sigma$. Note that
$\sS_{<\sigma}=\{\tau\in\sS\colon\tau<\sigma\}$ is the face poset
of the boundary of the simplex~$\Delta^\sigma$. Geometrically, we
may think of $|\sS|$ as obtained from $|\mathcal T|$ by attaching
one simplex $\Delta^\sigma$ along its boundary (if $|\sigma|=1$,
then $\Delta^\sigma$ is a single point, so $|\sS|$ is a disjoint
union of $|\mathcal T|$ and a point). We therefore need to prove
that the following is a pullback diagram:
\begin{equation}\label{stpb}
\begin{CD}
  \Z[\sS] @>>> \Z[\sigma]@.=\Z[\sS_{\le\sigma}]\\
  @VVV @VVV\\
  \Z[\mathcal T] @>>> \Z[\sS_{<\sigma}].
\end{CD}
\end{equation}
Here the vertical arrows are obtained by mapping $v_\sigma$ to 0,
while the horizontal ones are obtained by mapping $v_\tau$ to 0
for $\tau\not\le\sigma$. Denote the pullback of~\eqref{stpb} by
$A$; we need to show that $\Z[\sS]\to A$ is an isomorphism.

Since the limits in $\cgr$ are created in the underlying category
$\cgg$ of graded abelian groups (graded $\Z$-modules), the
underlying group of $A$ is the direct sum of $\Z[\mathcal T]$ and
$\Z[\sigma]$ with the pieces $\Z[\sS_{<\sigma}]$ identified in
both groups. In other words,
\begin{equation}\label{3dec}
  A=T\oplus\Z[\sS_{<\sigma}]\oplus S,
\end{equation}
where $T$ is the complement to $\Z[\sS_{<\sigma}]$ in $\Z[\mathcal
T]$, and $S$ is the complement to $\Z[\sS_{<\sigma}]$ in
$\Z[\sigma]$. By Lemma~\ref{chain}, the group $\Z[\sS_{<\sigma}]$
has basis of monomials
$v_{\tau_1}^{\alpha_1}v_{\tau_2}^{\alpha_2}\cdots
v_{\tau_k}^{\alpha_k}$ with $\tau_k<\sigma$. Similarly, $S$ has
basis of those monomials with $\tau_k=\sigma$ and $\alpha_k>0$,
while $T$ has basis of those monomials with $\tau_k\not\le\sigma$
and $\alpha_k>0$. Yet another application of Lemma~\ref{chain}
gives a decomposition of~$\Z[\sS]$ identical to~\eqref{3dec}: a
basis element $v_{\tau_1}^{\alpha_1}v_{\tau_2}^{\alpha_2}\cdots
v_{\tau_k}^{\alpha_k}$ with $\alpha_k>0$ has either
$\tau_k\not\le\sigma$, or $\tau_k<\sigma$, or $\tau_k=\sigma$.
These three possibilities map to $T$, $\Z[\sS_{<\sigma}]$ and $S$
respectively. It follows that $\Z[\sS]\to A$ is a group
isomorphism. Since it is a ring map, it is also a ring
isomorphism, thus finishing the proof.
\end{proof}
\end{construction}

The description of $\Z[\sS]$ as a limit has the following
corollary, describing the functorial properties of the face ring.

\begin{proposition}\label{funct}
Let $f\colon\sS\to\mathcal T$ be a rank-preserving map of
simplicial posets. Define a homomorphism
\[
  f^*\colon\Z[w_\tau\colon
  \tau\in\mathcal T]\to\Z[v_\sigma\colon
  \sigma\in\mathcal S]
\]
by $f^*(w_\tau)=\sum_{\sigma\in f^{-1}(\tau)}v_\sigma$. Then $f^*$
descends to a ring homomorphism $\Z[\mathcal T]\to\Z[\sS]$, which
we continue to denote by~$f^*$.
\end{proposition}
\begin{proof}
The poset map $f$ gives rise to a functor
$f\colon\cat^{op}(\sS)\to\cat^{op}(\mathcal T)$ and therefore to
\[
  f^*\colon[\cat^{op}(\mathcal T),\cgr]\to[\cat^{op}(\sS),\cgr],
\]
where $[\cat^{op}(\sS),\cgr]$ denotes the functors from
$\cat^{op}(\sS)$ to $\cgr$. It is easy to see that
$f^*\Z[\cdot]_\mathcal T=\Z[\cdot]_\mathcal S$ (see
Lemma~\ref{zslim}), so we have the induced map of limits
$f^*\colon\Z[\mathcal T]\to\Z[\sS]$. We also have that
$f^*(w_\tau)=\sum_{\sigma\in f^{-1}(\tau)}v_\sigma$ by the
construction of $\lim$ in $\cgr$.
\end{proof}

\begin{example}
The folding map~\eqref{fold} induces a monomorphism
$\Z[\sK_\sS]\to\Z[\sS]$, which embeds $\Z[\sK_\sS]$ in $\Z[\sS]$
as the subring generated by the elements~$v_i$.
\end{example}

\begin{remark}
The functoriality property for the face ring $\Z[\sK]$ of a
simplicial complex was observed in~\cite[Prop.~3.4]{bu-pa02}.
However an attempt to prove Proposition~\ref{funct} directly from
the definition, by showing that $f^*(\mathcal I_{\mathcal
T})\subset\mathcal I_\sS$, runs into a complicated combinatorial
analysis of the poset structure. This is an example of a situation
where the use of an abstract categorical description of $\Z[\sS]$
proves to be beneficial.

The $\lim$-construction of $\Z[\sS]$ also opens the way to further
generalisations of the face ring, to more general posets and maybe
to simplicial sets. Whether these rings would have a nice
algebraic description like that of Definition~\ref{deffr} is
questionable though.
\end{remark}

\section{Moment-angle complexes}
Let $D^2$ denote the standard unit 2-disc and $S^1$ its boundary
circle. We further consider the unit \emph{polydisc} $(D^2)^m$ in
the complex space $\C^m$:
\[
  (D^2)^m=\bigl\{ (z_1,\ldots,z_m)\in\C^m\colon |z_i|\le1,\quad i=1,\ldots,m
  \bigr\}.
\]
For every $\sigma\in\sS$, consider the following subset in
$(D^2)^m$:
\[
  B_\sigma=\bigl\{(z_1,\ldots,z_m)\in
  (D^2)^m\colon |z_j|=1\text{ if }j\not\le\sigma\bigl\},
\]
Then $B_\sigma$ is homeomorphic to a product of $|\sigma|$ discs
and $m-|\sigma|$ circles. We have an inclusion $B_\tau\subset
B_\sigma$ whenever $\tau\le\sigma$. It follows that the assignment
$\sigma\mapsto B_\sigma$ defines a diagram from $\cat(\sS)$ to
$\top$, which we denote $(D^2,S^1)^\sS$.

\begin{definition}\label{defmac}
The \emph{moment-angle complex} corresponding to a simplicial
poset $\sS$ is
\begin{equation}\label{zscolim}
  \zs=\colim(D^2,S^1)^\sS.
\end{equation}
\end{definition}

The space $\zs$ is glued from the ``moment-angle blocks''
$B_\sigma$ according to the poset relation in $\sS$. When $\sS$ is
a simplicial complex $\sK$ it becomes the standard moment-angle
complex $\zk$ of~\cite[\S6.2]{bu-pa02}.

\begin{remark}
The definition of $\zs$ is readily generalised to an arbitrary
pair of spaces $(X,W)$ as $\zs(X,W):=\colim(X,W)^\sS$. An easy
argument similar to~\cite[Prop.~3.5]{pano08l} shows that
\[
  H^*\bigl(\colim(\C P^\infty,pt)^\sS;\Z\bigr)\cong\Z[\sS].
\]
\end{remark}

\begin{example}
Let $\sS$ be the simplicial poset of Fig.~\ref{sccfig}~(a). Then
$\zs$ is obtained by gluing two copies of $D^2\times D^2$ along
their boundary $S^3=D^2\times S^1\cup S^1\times D^2$. Therefore,
$\zs\cong S^4$. Here, $\sK_\sS=\Delta^1$ (a segment), and the
moment-angle complex map induced by~\eqref{fold} folds $S^4$
onto~$D^4$. Similarly, if $\sS$ is of Fig.~\ref{sccfig}~(b), then
$\zs\cong S^6$. Note that even-dimensional spheres do not appear
as moment-angle complexes~$\zk$ for simplicial complexes~$\sK$.
\end{example}

The polydisc $(D^2)^m$ has the natural coordinatewise action of
the $m$-torus $T^m$, with quotient the $m$-cube $I^m$. Since every
inclusion $B_\tau\subset B_\sigma$ is $T^m$-equivariant, the
moment-angle complex $\zs$ acquires a $T^m$-action.

The \emph{join} of simplicial posets $\sS_1$ and $\sS_2$ is the
simplicial poset $\sS_1\mathbin{*}\sS_2$ whose elements are pairs
$(\sigma_1,\sigma_2)$, with
$(\sigma_1,\sigma_2)\le(\tau_1,\tau_2)$ whenever
$\sigma_1\le\tau_1$ in $\sS_1$ and $\sigma_2\le\tau_2$ in $\sS_2$.
The following properties of $\zs$ are similar to those of $\zk$
and can be proved in a very much similar fashion
(see~\cite[Ch.~6]{bu-pa02}).

\begin{proposition}
\begin{itemize}
\item[(a)] $\mathcal Z_{\sS_1\mathbin{*}\sS_2}\cong
\mathcal Z_{\sS_1}\times\mathcal Z_{\sS_2}$;

\item[(b)] the quotient $\zs/T^m$ is homeomorphic to the cone over
$|\sS|$;

\item[(c)] if $|\sS|\cong S^{n-1}$, then $\zs$ is a manifold of dimension~$m+n$.
\end{itemize}
\end{proposition}

An important series of examples of simplicial posets $\sS$ with
$|\sS|\cong S^{n-1}$ comes from the inverse face posets of
\emph{face-acyclic manifolds with corners} in the sense
of~\cite{ma-pa06}. These manifolds with corners $Q$ provide
decompositions of an $n$-dimensional ball into faces, generalising
those face decompositions coming from \emph{simple
$n$-polytopes}~$P$. We therefore obtain \emph{moment-angle
mani\-folds} $\mathcal Z_Q$ generalising the manifolds $\mathcal
Z_P$ corresponding to simple polytopes.

\begin{construction}[cell decomposition]
The disc $D^2$ decomposes in the standard way into three cells of
dimensions $0$, $1$ and~$2$, which we denote $*$, $T$ and $D$
respectively. The polydisc $(D^2)^m$ then acquires the product
cell decomposition, with each $B_\tau\subset B_\sigma$ being an
inclusion of cellular subcomplexes for $\tau\le\sigma$. We
therefore obtain a cell decomposition of $\zs$. Each cell in $\zs$
is determined by an element $\sigma\in\sS$ and a subset $\omega\in
V(\sS)$ with $V(\sigma)\cap\omega=\varnothing$. Such a cell is a
product of $|\sigma|$ cells of $D$-type, $|\omega|$ cells of
$T$-type and the rest of $*$-type. We denote this cell by
$\kappa(\omega,\sigma)$.

The resulting cellular cochain complex $C^*(\zs)$ has an additive
basis consisting of cochains $\kappa(\omega,\sigma)^*$ dual to the
corresponding cells. We introduce a $\Z\oplus\Z^m$-grading on the
cochains by setting
\[
  \mdeg\kappa(\omega,\sigma)^*=(-|\omega|,2V(\sigma)+2\omega),
\]
where we think of both $V(\sigma)$ and $\omega$ as vectors
in~$\{0,1\}^m\subset\Z^m$. The cellular differential does not
change the $\Z^m$-part of the multigrading, so we obtain a
decomposition
\[
  C^*(\zs)=\bigoplus_{\mb a\in\Z^m} C^{*,2\mb a}(\zk)
\]
into a sum of subcomplexes. In fact the only nontrivial
subcomplexes are those for which $\mb a$ is in~$\{0,1\}^m$. The
cellular cohomology of $\zs$ thereby acquires an additional
grading, and we may define the \emph{multigraded Betti numbers}
$b^{-i,2\mb a}(\zs)$ by
\[
  b^{-i,2\mb a}(\zs)=\mathop{\mathrm{rank}} H^{-i,2\mb a}(\zs),
  \quad i=1,\ldots,m,\;\mb a\in\Z^m.
\]
For the ordinary Betti numbers we have $b^k(\zs)=\sum_{2|\mb
a|-i=k}b^{-i,2\mb a}(\zs)$.
\end{construction}

The face ring $\Z[\sS]$ acquires a $\Z[v_1,\ldots,v_m]$-algebra
structure via the map $\Z[v_1,\ldots,v_m]\to\Z[\sS]$ sending each
$v_i$ identically. (Unlike the case of simplicial complexes, this
map is generally not surjective.) The $\Z\oplus\Z^m$-graded
$\Tor$-algebra of $\Z[\sS]$ is defined in the standard
way~\cite[\S3.4]{bu-pa02}:
\[
  \Tor_{\Z[v_1,\ldots,v_m]}(\Z[\sS],\Z)=
  \bigoplus_{i\ge0,\mb a\in\Z^m}
  \Tor^{-i,\,2\mb a}_{\Z[v_1,\ldots,v_m]}(\Z[\sS],\Z).
\]
Note that the first degree is always nonpositive, which is because
we number the terms in a free resolution by nonpositive integers.

\begin{theorem}\label{hzs1}
There is a graded ring isomorphism
\[
  H^*(\zs;\Z)\cong\Tor_{\Z[v_1,\ldots,v_m]}(\Z[\sS],\Z)
\]
whose graded components are given by the group isomorphisms
\begin{equation}\label{additive}
  H^p(\zs;\Z)\cong\bigoplus_{-i+2|\mb a|=p}
  \Tor^{-i,\,2\mb a}_{\Z[v_1,\ldots,v_m]}(\Z[\sS],\Z)
\end{equation}
in each degree $p$. Here $|\mb a|=j_1+\cdots+j_m$ for $\mb
a=(j_1,\ldots,j_m)$.
\end{theorem}

Using the Koszul resolution
for the trivial $\Z[v_1,\ldots,v_m]$-module $\Z$ (see
\cite[Lemma~3.29]{bu-pa02}) we restate the above theorem as
follows:

\begin{theorem}\label{hzs2}
There is a graded ring isomorphism
\[
  H^*(\zs;\Z)\cong
  H\bigl[\Lambda[u_1,\ldots,u_m]\otimes\Z[\sS],d\bigr].
\]
Here on the right hand side stands the cohomology of a
differential $\Z\oplus\Z^m$-graded ring with
\[
  \mdeg u_i=(-1,2\mb e_i),\quad \mdeg v_\sigma=(0,2V(\sigma)),
  \quad du_i=v_i, \quad dv_\sigma=0,
\]
where $\mb e_i\in\Z^m$ it the $i$th basis vector, for
$i=1,\ldots,m$.
\end{theorem}
\begin{proof}
The proof  given here structurally resembles the proof
of~\cite[Th.~4.7]{pano08l} (for the case of~$\zk$). However,
algebraic arguments used in the proof for $\zk$ do not work in the
case of simplicial posets. Instead, we use topological and
categorical arguments in the appropriated places of this proof.

We consider the quotient differential graded ring
\[
  R^*(\sS):=\Lambda[u_1,\ldots,u_m]\otimes\Z[\sS]/\mathcal I_R
\]
where $\mathcal I_R$ is the ideal generated by the elements
\[
  u_iv_\sigma\quad\text{with }i\in V(\sigma),\quad\text{and}\quad
  v_\sigma v_\tau\;\;\;\text{with }\sigma\wedge\tau\ne\hatzero.
\]
Note that the latter condition is equivalent to $V(\sigma)\cap
V(\tau)\ne\varnothing$.

We claim that the quotient projection
\[
  \varrho\colon\Lambda[u_1,\ldots,u_m]\otimes\Z[\sS]\to R^*(\sS)
\]
is a quasi-isomorphism, that is, it induces an isomorphism in
cohomology.

Lemma~\ref{chain} implies that $R^*(\sS)$ is generated, as an
abelian group, by the monomials $u_\omega v_\sigma$ where
$\omega\subseteq V(\sS)$, \ $\sigma\in\sS$, \ $\omega\cap
V(\sigma)=\varnothing$, and $u_\omega=u_{i_1}\ldots u_{i_k}$ for
$\omega=\{i_1,\ldots,i_k\}$. In particular, $R^*(\sS)$ is a free
abelian group of finite rank. It is now easy to observe that the
map
\begin{align}\label{gmap}
  g\colon R^*(\sS) &\to C^*(\zs),\\
  u_\omega v_\sigma &\mapsto\mathcal \kappa(\omega,\sigma)^*\notag
\end{align}
is an isomorphism of cochain complexes. Indeed, the additive bases
of the two groups are in one-to-one correspondence, and the
differential in $R^*(\sS)$ acts (in the case $|\omega|=1$ and
$i\notin V(\sigma)$) as
\[
  d(u_iv_\sigma)=v_iv_\sigma=\sum_{\eta\in i\vee\sigma}v_\eta.
\]
This is exactly how the cellular differential in $C^*(\zs)$ acts
on $\kappa(i,\sigma)^*$. The case of arbitrary $\omega$ is treated
similarly. It follows that we have an isomorphism of cohomology
groups $H^j[R^*(\sS)]\cong H^j(\zs)$ for all~$j$.

The differential ring $\Lambda[u_1,\ldots,u_m]\otimes\Z[\sS]$ also
may be identified with the cellular cochains of a certain space.
Namely, consider the space $\zs(S^\infty,S^1)$ defined in the same
way as~\eqref{zscolim}, but with $D^2$ replaced by an
infinite-dimensional sphere~$S^\infty$. The latter is a
contractible space which has a cell decomposition with one cell in
every dimension. The boundary of every $2k$-dimensional cell is
the closure of the $(2k-1)$-cell, while the boundary of an
odd-dimensional cell is zero. The cellular cochains of $S^\infty$
can be identified with the Koszul differential ring
\[
  \Lambda[u]\otimes\Z[v],\quad\deg u=1,\;\deg v=2,\quad du=v,\;dv=0.
\]
As in the case of~\eqref{gmap}, Lemma~\ref{chain} implies that
there is an isomorphism of cochain complexes
\[
  g'\colon\Lambda[u_1,\ldots,u_m]\otimes\Z[\sS]\to C^*(\zs(S^\infty,S^1)).
\]
We also have a deformation retraction $D^2\hookrightarrow
S^\infty\to D^2$. It follows from the standard functoriality
arguments that we also have a deformation retraction
\[
  \zs=\colim(D^2,S^1)^\sS\hookrightarrow\colim(S^\infty,S^1)^\sS
  \to\colim(D^2,S^1)^\sS
\]
onto a cellular subcomplex. Therefore the cochain map
$C^*(\zs(S^\infty,S^1))\to C^*(\zs)$ induced by the inclusion is a
cohomology isomorphism.

Summarising the above observations we obtain the commutative
square
\begin{equation}\label{comsq}
\begin{CD}
  \Lambda[u_1,\ldots,u_m]\otimes\Z[\sS] @>g'>>
  C^*(\zs(S^\infty,S^1))\\
  @V\varrho VV @VVV\\
  R^*(\sS) @>g>> C^*(\zs)
\end{CD}
\end{equation}
in which the horizontal arrows are isomorphisms of cochain
complexes, and the right vertical arrow induces a cohomology
isomorphism. It follows that the left arrow is a
quasi-isomorphism, as claimed.

\begin{remark}
There is an obvious inclusion of cochain complexes $\iota\colon
R^*(\sS)\to\Lambda[u_1,\ldots,u_m]\otimes\Z[\sS]$, which is not a
ring homomorphism though. It is possible to prove that $\varrho$
is a cohomology isomorphism by constructing a cochain homotopy $s$
between the maps $\mathop\mathrm{id}$ and $\iota\cdot\varrho$ from
$\Lambda[u_1,\ldots,u_m]\otimes\Z[\sS]$ to itself. However, in the
construction of $s$ we cannot use an inductive argument as
in~\cite[Lemma~4.4]{pano08l}, and the general formula for $s$ is
rather cumbersome.
\end{remark}

The additive isomorphism of~\eqref{additive} now follows
from~\eqref{comsq}. To establish the ring isomorphism we need to
analyse the multiplication of cellular cochains in~$C^*(\zs)$.

We consider the diagonal approximation map
$\widetilde{\Delta}\colon D^2\to D^2\times D^2$, defined in polar
coordinates $z=\rho e^{i\varphi}\in D^2$, $0\le\rho\le1$,
$0\le\varphi<2\pi$, as follows:
\[
  \rho e^{i\varphi}\mapsto\left\{
  \begin{array}{ll}
    (1+\rho(e^{2i\varphi}-1),1)&\text{ for }0\le\varphi\le\pi,\\
    (1,1+\rho(e^{2i\varphi}-1))&\text{ for }\pi\le\varphi<2\pi.
  \end{array}
  \right.
\]
This is a cellular map homotopic to the diagonal $\Delta\colon
D^2\to D^2\times D^2$. Taking an $m$-fold product, we obtain a
cellular diagonal approximation
\[
  \widetilde\Delta\colon(D^2)^m\to(D^2)^m\times(D^2)^m.
\]
It restricts to a map $B_\sigma\to B_\sigma\times B_\sigma$ for
every $\sigma\in\sS$ and gives rise to a map of diagrams
\[
  (D^2,S^1)^\sS\to(D^2,S^1)^\sS\times(D^2,S^1)^\sS.
\]
By definition, the colimit of the latter is $\mathcal
Z_{\sS\mathop{*}\sS}$, which is identified with $\zs\times\zs$. We
therefore obtain a cellular approximation
$\widetilde\Delta\colon\zs\to\zs\times\zs$ for the diagonal map
of~$\zs$. It induces a ring structure on the cellular cochains via
the composition
\[
\begin{CD}
  C^*(\zs)\otimes C^*(\zs) @>\times>> C^*(\zs\times\zs)
  @>\widetilde{\Delta}^*>> C^*(\zs).
\end{CD}
\]
We claim that, with this multiplication in $C^*(\zs)$, the
map~\eqref{gmap} becomes a differential graded ring isomorphism.
To see this we first observe that since~\eqref{gmap} is a linear
map, it is enough to consider the product of two generators
$u_\omega v_\sigma$ and $u_\psi v_\tau$. If any two of the subsets
$\omega$, $V(\sigma)$, $\psi$ and $V(\tau)$ have nonempty
intersection, then $u_\omega v_\sigma\cdot u_\psi v_\tau=0$.
Otherwise (if all of the four subsets are complementary) we have
\begin{equation}\label{gprod}
  g(u_\omega v_\sigma\cdot u_\psi v_\tau)=
  g\bigl(u_{\omega\sqcup\psi}\cdot
  \sum_{\eta\in\,\sigma\vee\tau}v_\eta\bigr)=
  \sum_{\eta\in\,\sigma\vee\tau}\kappa(\omega\sqcup\psi,\eta)^*.
\end{equation}
We also observe that
\[
  \widetilde\Delta\kappa(\chi,\eta)=
  \mathop{\sum_{\omega\sqcup\psi=\chi}}\limits_{\sigma\vee\tau\,\ni\;\eta}
  \kappa(\omega,\sigma)\times\kappa(\psi,\tau)
\]
whenever $\chi\cap V(\eta)=\varnothing$. Therefore,
\begin{multline*}
  g(u_\omega v_\sigma)\cdot g(u_\psi v_\tau)=
  \kappa(\omega,\sigma)^*\cdot\kappa(\psi,\tau)^*\\=
  \widetilde\Delta^*
  \bigl(\kappa(\omega,\sigma)\times\kappa(\psi,\tau)\bigr)^*=
  \sum_{\eta\,\in\,\sigma\vee\tau}\kappa(\omega\sqcup\psi,\eta)^*.
\end{multline*}
Comparing this with~\eqref{gprod} we deduce that~\eqref{gmap}
is a ring map, concluding the proof.
\end{proof}

\begin{remark}
Using the monoid structure on $D^2$ as
in~\cite[Lemma~4.2]{pano08l} one easily sees that the construction
of $\zs$ is functorial with respect to maps of simplicial posets.
This together with Proposition~\ref{funct} makes the isomorphism
of Theorem~\ref{hzs1} functorial.
\end{remark}

We have the following important corollary.

\begin{corollary}
The groups $\Tor^{-i,\,2\mb a}_{\Z[v_1,\ldots,v_m]}(\Z[\sS],\Z)$
vanish for $\mb a\notin\{0,1\}^m$.
\end{corollary}
\begin{proof}
The multigraded component $R^{-i,\,2\mb a}(\sS)$ is zero for $\mb
a\notin\{0,1\}^m$.
\end{proof}

Denote by $\sS_{\mb a}$ the
subposet of $\sS$ consisting of those $\sigma$ for which
$V(\sigma)\subset\mb a$. As a further corollary of
Theorems~\ref{hzs1} and~\ref{hzs2} we obtain the following
generalisation of Hochster's theorem to simplicial posets.

\begin{corollary}\label{gh}
For every $\mb a\in\{0,1\}^m$ there is an isomorphism
\[
  \Tor^{-i,\,2\mb a}_{\Z[v_1,\ldots,v_m]}(\Z[\sS],\Z)\cong
  \widetilde H^{|\mb a|-i-1}\bigl(|\sS_{\mb a}|\bigr)
\]
(here we follow the standard convention that $\widetilde
H^{-1}(\varnothing)=\Z$).
\end{corollary}
\begin{proof}
The argument is identical to that of~\cite[Th.~5.1]{pano08l}:
there is an isomorphism of cellular cochain complexes
\[
  \widetilde C^*\bigl(|\sS_{\mb a}|\bigr)\to
  C^{*\,+1-|\mb a|,\,2\mb a}(\zk),
  \quad \sigma^*\mapsto
  \kappa\bigl(\mb a\!\setminus\! V(\sigma),\,\sigma\bigr)^*,
\]
inducing the required isomorphisms in cohomology.
\end{proof}

\begin{remark}
The statement of Corollary~\ref{gh} was obtained by
Duval~\cite{duva97} (with field coefficients, and without
considering the ring structure in $\Tor$).

It is clear from Corollary~\ref{gh} that the cohomology of $\zs$
may contain an arbitrary amount of additive torsion; just take
$|\sS|$ to be a triangulation of a space with the appropriate
torsion in cohomology.
\end{remark}

The \emph{multigraded algebraic Betti numbers} of $\Z[\sS]$ are
defined as
\[
  \beta^{-i,\,2\mb a}(\Z[\sS])=
  \mathop\mathrm{rk}\Tor^{-i,\,2\mb a}_{\Z[v_1,\ldots,v_m]}(\Z[\sS],\Z)
  =\mathop\mathrm{rk}H^{-i,\,2\mb a}(\zs)
\]
for $i=1,\ldots,m$, $\mb a\in\Z^m$. We also set
$\beta^{-i}(\Z[\sS])=\sum_{\mb a\in\Z^m}\beta^{-i,\,2\mb
a}(\Z[\sS])$.

\begin{example}
Let us see how the isomorphism of Theorem~\ref{hzs2} looks in the
case of the simplicial poset $\sS$ of Example~\ref{exsp}.1. The
elements $1$, $v_1$, $v_2$, $v_\sigma$ and $v_\tau$ of $R^{0,*}$
are all cocycles. Moreover, $v_1$, $v_2$ and $v_\sigma+v_\tau$ are
coboundaries, the latter because
$d(u_1v_2)=v_1v_2=v_\sigma+v_\tau$. It therefore follows that
$\beta^{0,(0,0)}(\Z[\sS])=\beta^{0,(2,2)}(\Z[\sS])=1$, while
$\beta^{0,(2,0)}(\Z[\sS])=\beta^{0,(0,2)}(\Z[\sS])=0$. Also, a
direct computation shows that $\beta^{-i,2\mb a}(\Z[\sS])=0$ for
$i>0$. This implies that $\Z[\sS]$ is a free $\Z[v_1,v_2]$-module
with two generators, $1$ and $v_\sigma$. The multigraded
decomposition~\eqref{additive} in cohomology of $\zs\cong S^4$ is
as follows: $H^{0}(\zs)=H^{0,(0,0)}(\zs)\cong\Z$ and
$H^{4}(\zs)=H^{0,(2,2)}(\zs)\cong\Z$.

The reader may compare this with similar computations
of~\cite[Ex.~4.8, Ex.~5.7]{pano08l} in the case of moment-angle
complexes $\zk$. Note that unlike the case of simplicial
complexes, $\beta^0(\Z[\sS])$ may be bigger than 1. In fact,
Corollary~\ref{gh} implies the following.
\begin{proposition}
The number of generators of $\Z[\sS]$ as a
$\Z[v_1,\ldots,v_m]$-module equals
\[
  \beta^0(\Z[\sS])=\sum_{\mb a\subset[m]}
  \rk\widetilde H^{|\mb a|-1}\bigl(|\sS_{\mb a}|\bigr).
\]
\end{proposition}
\end{example}

We finish this section by considering a poset $\sS$ slightly more
complicated than the toy examples we saw before, and calculating
the cohomology of $\zs$ accordingly.

\begin{example}\label{exmwc}
\begin{figure}
\begin{picture}(0,0)%
\includegraphics{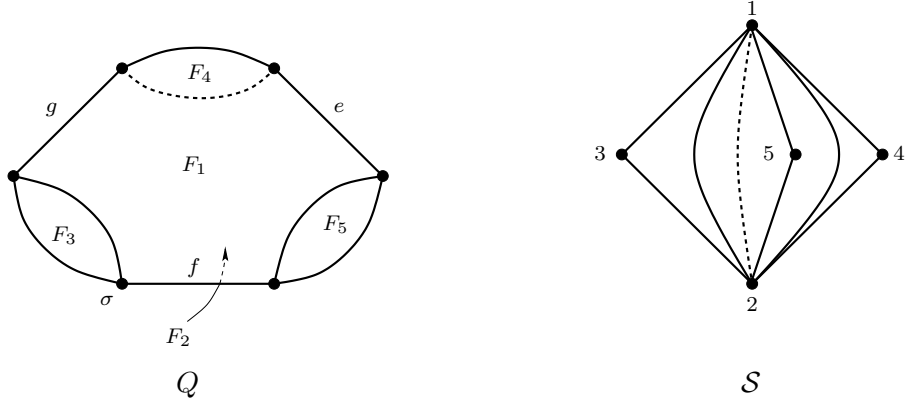}%
\end{picture}%
\setlength{\unitlength}{1776sp}%
\begingroup\makeatletter\ifx\SetFigFont\undefined%
\gdef\SetFigFont#1#2#3#4#5{%
  \reset@font\fontsize{#1}{#2pt}%
  \fontfamily{#3}\fontseries{#4}\fontshape{#5}%
  \selectfont}%
\fi\endgroup%
\begin{picture}(12648,5550)(8018,-6769)
\put(18151,-6661){\makebox(0,0)[lb]{\smash{{\SetFigFont{11}{13.2}{\rmdefault}{\mddefault}{\updefault}$\mathcal{S}$}}}}
\put(10201,-5986){\makebox(0,0)[lb]{\smash{{\SetFigFont{8}{9.6}{\rmdefault}{\mddefault}{\updefault}$F_2$}}}}
\put(12526,-2761){\makebox(0,0)[lb]{\smash{{\SetFigFont{8}{9.6}{\rmdefault}{\mddefault}{\updefault}$e$}}}}
\put(10501,-2311){\makebox(0,0)[lb]{\smash{{\SetFigFont{8}{9.6}{\rmdefault}{\mddefault}{\updefault}$F_4$}}}}
\put(8626,-4561){\makebox(0,0)[lb]{\smash{{\SetFigFont{8}{9.6}{\rmdefault}{\mddefault}{\updefault}$F_3$}}}}
\put(12376,-4411){\makebox(0,0)[lb]{\smash{{\SetFigFont{8}{9.6}{\rmdefault}{\mddefault}{\updefault}$F_5$}}}}
\put(20251,-3436){\makebox(0,0)[lb]{\smash{{\SetFigFont{8}{9.6}{\rmdefault}{\mddefault}{\updefault}$4$}}}}
\put(18226,-5536){\makebox(0,0)[lb]{\smash{{\SetFigFont{8}{9.6}{\rmdefault}{\mddefault}{\updefault}$2$}}}}
\put(18451,-3436){\makebox(0,0)[lb]{\smash{{\SetFigFont{8}{9.6}{\rmdefault}{\mddefault}{\updefault}$5$}}}}
\put(18226,-1411){\makebox(0,0)[lb]{\smash{{\SetFigFont{8}{9.6}{\rmdefault}{\mddefault}{\updefault}$1$}}}}
\put(8551,-2761){\makebox(0,0)[lb]{\smash{{\SetFigFont{8}{9.6}{\rmdefault}{\mddefault}{\updefault}$g$}}}}
\put(16126,-3436){\makebox(0,0)[lb]{\smash{{\SetFigFont{8}{9.6}{\rmdefault}{\mddefault}{\updefault}$3$}}}}
\put(10426,-3586){\makebox(0,0)[lb]{\smash{{\SetFigFont{8}{9.6}{\rmdefault}{\mddefault}{\updefault}$F_1$}}}}
\put(9301,-5461){\makebox(0,0)[lb]{\smash{{\SetFigFont{8}{9.6}{\rmdefault}{\mddefault}{\updefault}$\sigma$}}}}
\put(10501,-5011){\makebox(0,0)[lb]{\smash{{\SetFigFont{8}{9.6}{\rmdefault}{\mddefault}{\updefault}$f$}}}}
\put(10351,-6661){\makebox(0,0)[lb]{\smash{{\SetFigFont{11}{13.2}{\rmdefault}{\mddefault}{\updefault}$Q$}}}}
\end{picture}%
\centering
  \caption{Manifold with corners $Q$ and the dual poset~$\sS$.}
\label{mwc}
\end{figure}
Let $Q$ be the 3-dimensional manifold with corners shown in
Figure~\ref{mwc}~(left). It is a 3-ball with $m=5$ facets
$F_1,\ldots,F_5$ numbered as shown. We denote the edges $e,f,g$
and the vertex $\sigma$ of $Q$ as shown. The corresponding
moment-angle complex $\zq$ is an 8-dimensional manifold.

The inverse face poset of $Q$ is the simplicial poset~$\sS$ shown
in Figure~\ref{mwc} (right). Note that the facets of $Q$
correspond to the 5 vertices of~$\sS$, while $\sigma$ corresponds
to a certain 2-simplex of~$\sS$. The face ring $\Z[\sS]$ is the
quotient of the polynomial ring
\[
  \Z[\sS]=\Z[v_1,\ldots,v_5,v_e,v_f,v_g],\quad\deg v_i=2,\quad
  \deg v_e=\deg v_f=\deg v_e=4
\]
by the relations
\begin{align*}
  v_1v_2&=v_e+v_f+v_g,\\
  v_3v_4&=v_3v_5=v_4v_5=
  v_3v_e=v_4v_f=v_5v_g=v_ev_f=v_ev_g=v_ev_f=0.
\end{align*}
The other generators and relations in the original presentation
can be derived from these; e.g., $v_\sigma=v_3v_f$.

Given a vector $\mb a\in\{0,1\}^m$ regarded as a subset of~$[m]$,
set
\[
  Q_\mb a=\bigcup_{i\in\mb a}F_i\subset Q.
\]
It is a subspace in the boundary of~$Q$. Using the barycentric
subdivision it is easy to see that $|\sS_{\mb a}|$ is a
deformation retract of~$Q_\mb a$. Then Theorem~\ref{hzs1} and
Corollary~\ref{gh} give the following formula for the multigraded
cohomology of~$\zq$:
\begin{equation}\label{faces}
  H^{-i,\,2\mb a}(\zq)\cong
  \widetilde H^{|\mb a|-i-1}(Q_\mb a).
\end{equation}
Using this formula we calculate the nontrivial cohomology groups
of $\zq$ as follows:
\begin{align*}
H^{0,(0,0,0,0,0)}(\zq)&=\widetilde H^{-1}(\varnothing)=\Z
 &&1\\
H^{-1,(0,0,2,2,0)}(\zq)&=\widetilde H^{0}(F_3\cup F_4)=\Z
 &&u_3v_4\\
H^{-1,(0,0,2,0,2)}(\zq)&=\widetilde H^{0}(F_3\cup F_5)=\Z
 &&u_5v_3\\
H^{-1,(0,0,0,2,2)}(\zq)&=\widetilde H^{0}(F_4\cup F_5)=\Z
 &&u_4v_5\\
H^{-2,(0,0,2,2,2)}(\zq)&=\widetilde H^{0}(F_3\cup F_4\cup
F_5)=\Z\oplus\Z
 &&u_5u_3v_4,\;u_5u_4v_3\\
H^{0,(2,2,0,0,0)}(\zq)&=\widetilde H^{1}(F_1\cup F_2)=\Z\oplus\Z
 &&v_e,\;v_f\\
H^{-1,(2,2,2,0,0)}(\zq)&=\widetilde H^{1}(F_1\cup F_2\cup F_3)=\Z
 &&u_3v_e\\
H^{-1,(2,2,0,2,0)}(\zq)&=\widetilde H^{1}(F_1\cup F_2\cup F_4)=\Z
 &&u_4v_f\\
H^{-1,(2,2,0,0,2)}(\zq)&=\widetilde H^{1}(F_1\cup F_2\cup F_5)=\Z
 &&u_5v_g\\
H^{-2,(2,2,2,2,2)}(\zq)&=\widetilde H^{2}(F_1\cup\cdots\cup
F_5)=\Z && u_5u_4v_3v_f=u_5u_4v_\sigma
\end{align*}
It follows that the ordinary (1-graded) Betti numbers of $\zq$ are
given by the sequence $(1,0,0,3,4,3,0,0,1)$. In the right column
of the table above we include the cocycles in the differential
graded ring $\Lambda[u_1,\ldots,u_5]\otimes\Z[\sS]$ representing
generators of the corresponding cohomology group. This allows us
to determine the ring structure in $H^*(\zs)$. For example,
\[
  [u_5u_3v_4]\cdot[v_f]=[u_5u_3v_4v_f]=0=[u_5u_4v_3]\cdot[v_e].
\]
On the other hand,
\begin{multline*}
  [u_5u_3v_4]\cdot[v_e]=-[u_3u_5v_4v_e]=-[u_3u_4v_5v_e]
  \\=[u_3u_4v_5v_f]=[u_5u_4v_3v_f]=[u_5u_4v_3]\cdot[v_f].
\end{multline*}
Here we have used the relations
$d(u_3u_4u_5v_e)=u_3u_4v_5v_e-u_3u_5v_4v_e$ and
$d(u_1u_3u_4v_2v_5)=u_3u_4v_5v_e+u_3u_4v_5v_f$. In fact, all
nontrivial products come from Poincar\'e duality. These
calculations may be summarised by the cohomology ring isomorphism
\[
  H^*(\zq)\cong H^*\bigl((S^3\times S^5)^{\#3}\cs
  (S^4\times S^4)^{\#2}\bigr)
\]
where the manifold on the right hand side is the connected sum of
three copies of $S^3\times S^5$ and two copies of $S^4\times S^4$.
We expect that this cohomology isomorphism is induced by a
homeomorphism; one might be able to prove this by using the
surgery techniques of~\cite{gi-lo09}.
\end{example}

\section{Almost free torus actions}
Halperin's \emph{toral rank conjecture} states that if a torus
$T^k$ acts almost freely on a finite-dimensional space $X$, then
the ``total amount of homology'' of $X$ is at least that of the
torus, that is,
\[
  \sum_i\rk H^i(X)\ge 2^k.
\]
(An action is \emph{almost free} if all isotropy subgroups are
finite.) We refer to $\sum_i\rk H^i(X)$ as the \emph{homology rank
of $X$} and denote it $\hrk(X)$.

It has been shown in the recent works of Cao--L\"u~\cite{ca-lu09}
and Ustinovsky~\cite{usti09} that the toral rank conjecture holds
for the restricted torus action on the moment-angle complex~$\zk$.
Here we show that the same holds for~$\zs$.

We define the \emph{toral rank} $\trk\zs$ as the maximal dimension
of a subtorus $T^k\subset T^m$ acting almost freely on~$\zs$.
Assume that $\dim\sS=n-1$; then $\dim\zs=m+n$. The isotropy
subgroups of the $T^m$-action on $\zs$ are coordinate subtori in
$T^m$ of the form
\begin{equation}\label{tsigma}
  T^\sigma=\{(z_1,\ldots,z_m)\in T^m\colon z_i=1\text{ for }
  i\notin V(\sigma)\}
\end{equation}
where $\sigma\in\sS$. The maximal dimension of these subgroups is
$n$, hence $\trk\zs\le m-n$.

Let $\mb t=(t_1,\dots,t_n)$ be a sequence of linear (degree-two)
elements in $\Z[\sS]$. We may write
\begin{equation}\label{lsop}
  t_i=\lambda_{i1}v_1+\cdots+\lambda_{im}v_m,\quad i=1,\ldots,n.
\end{equation}
Given $\sigma\in\sS$, define the \emph{restriction homomorphism}
\[
  s_\sigma\colon\Z[\sS]\to
  \Z[\sS]/(v_\tau\colon\tau\not\le\sigma).
\]
Its image may be identified with the polynomial ring $\Z[\sigma]$
on $|\sigma|$ generators. Note that $s_\sigma$ is induced by the
inclusion of posets $\sS_{\le\sigma}\to\sS$. Remember that $\mb t$
is called an \emph{lsop} (\emph{linear system of parameters})
in~$\Z[\sS]$ if it consists of algebraically independent elements
and $\Z[\sS]$ is a finitely generated $\Z[t_1,\ldots,t_n]$-module
(equivalently, $\Z[\sS]/(\mb t)$ has finite rank as an abelian
group).

\begin{lemma}\label{restr}
A degree-two sequence $\mb t=(t_1,\ldots,t_n)$ is an lsop in
$\Z[\sS]$ if and only if for every $\sigma\in\sS$ the elements
$s_\sigma(t_1),\ldots,s_\sigma(t_n)$ generate the positive degree
ideal $\Z[\sigma]_+$.
\end{lemma}
\begin{proof}
Assume \eqref{lsop} is an lsop. Every $s_\sigma$ induces an
epimorphism of the quotient rings:
\[
  \Z[\sS]/(\mb t)\to\Z[\sigma]/s_\sigma(\mb t).
\]
Since $\mb t$ is an lsop, $\Z[\sS]/(\mb t)$ has finite rank as a
group. Therefore, $\Z[\sigma]/s_\sigma(\mb t)$ is also of finite
rank, which happens only if $s_\sigma(\mb t)$ generates
$\Z[\sigma]_+$.

The other direction is proved by considering the sum of the
restrictions:
\[
  \Z[\sS]\longrightarrow\bigoplus_{\sigma\in\sS}\Z[\sigma].
\]
This is an injective $\Z[t_1,\ldots,t_n]$-module map
by~\cite[Lemma~5.6]{ma-pa06}. Since $\Z[t_1,\ldots,t_n]$ is a
Noetherian ring and $\bigoplus_{\sigma\in\sS}\Z[\sigma]$ is
finitely generated as a $\Z[t_1,\ldots,t_n]$-module by assumption,
its submodule $\Z[\sS]$ is also finitely generated. This implies
that $\mb t$ is an lsop.
\end{proof}


We organise the coefficients in \eqref{lsop} into an $n\times
m$-matrix $\Lambda=(\lambda_{ij})$. For any $\sigma\in\sS$ denote
by $\Lambda_\sigma$ the $n\times|\sigma|$-submatrix formed by the
elements $\lambda_{ij}$ with $j\in V(\sigma)$. The matrix
$\Lambda$ defines homomorphisms $\Z^m\to\Z^n$ and $\lambda\colon
T^m\to T^n$. Let $T_\Lambda=\ker\lambda\subset T^m$.

\begin{theorem}\label{zkquoq}
The following conditions are equivalent:
\begin{itemize}
\item[(a)] the sequence \eqref{lsop} is an lsop in the rational face
  ring $\Q[\sS]$;
\item[(b)] for every $\sigma\in\sS$ the matrix $\Lambda_\sigma$ has rank $|\sigma|$;
\item[(c)] $T_\Lambda$ is the product of an $(m-n)$-torus and a finite group,
 and $T_\Lambda$ acts almost freely on~$\zs$.
\end{itemize}
\end{theorem}
\begin{proof}
The equivalence of (a) and (b) is the $\Q$-version of
Lemma~\ref{restr}. Now, (b) holds if and only if $T_\Lambda\cap
T^\sigma$ is a finite group for every $\sigma\in\sS$, which means
that $T_\Lambda$ acts almost freely on~$\zs$ (see~\eqref{tsigma}).
The fact that $T_\Lambda$ contains an $(m-n)$-torus also follows
from~(b), because there is $\sigma\in\sS$ with $|\sigma|=n$.
\end{proof}

\begin{corollary}\label{trkzs}
If $\sS$ is of rank $n$ with $m$ vertices, then $\trk\zs=m-n$.
\end{corollary}
\begin{proof}
Consider the ring $\Q[\mathcal K_\mathcal S]$. Since it is
generated by the degree-two elements, it has an lsop $\mb t$ (this
is where we need the $\Q$-coefficients). Since $\Q[\sS]$ is
integral over $\Q[\mathcal K_\mathcal S]$
by~\cite[Lemma~3.9]{stan91}, $\mb t$ is also an lsop for
$\Q[\sS]$. By multiplying by a common denominator, we may assume
that $\mb t$ is in~$\Z[\sS]$ (although it may fail to be an
\emph{integral} lsop). Then there is an $(m-n)$-subtorus acting
almost freely on $\zs$ by Theorem~\ref{zkquoq}.
\end{proof}

There is also an integral version of Theorem~\ref{zkquoq}, which
is proved similarly:

\begin{theorem}
The following conditions are equivalent:
\begin{itemize}\label{quot}
\item[(a)] the sequence \eqref{lsop} is an lsop in $\Z[\sS]$;
\item[(b)] for every $\sigma\in\sS$ the columns of
$\Lambda_\sigma$ form a part of a basis of $\Z^n$;
\item[(c)] $T_\Lambda$ is an $(m-n)$-torus acting freely on~$\zs$.
\end{itemize}
\end{theorem}

\begin{remark}
Unlike the case of $\Q[\sS]$, an lsop in $\Z[\sS]$ may fail to
exist, which means the there is no $(m-n)$-subtorus acting
\emph{freely} on the corresponding~$\zs$. The maximal dimension
$s(\sS)$ of a subtorus $T^s\subset T^m$ acting freely on~$\zs$ is
also known as the \emph{Buchstaber invariant} of~$\sS$. It is a
much more subtle characteristic than $\trk(\zs)$ and is usually
difficult to determine. For more information about the Buchstaber
invariant for polytopes and simplicial complexes see~\cite{erok09}
and~\cite{fu-ma09}.
\end{remark}

\begin{proposition}
We have that $\hrk\zs\ge\hrk\mathcal Z_{\mathcal K_\mathcal S}$.
\end{proposition}
\begin{proof}
The folding map $|\sS|\to|\mathcal K_\mathcal S|$ has an obvious
section, which means that it is a retraction. It follows that $\rk
\widetilde H^i(|\sS|)\ge\rk\widetilde H^i(|\mathcal K_\mathcal
S|)$. The same holds for every subposet $\mathcal S_\mb a$. Now
the result follows from Theorem~\ref{hzs1} and Corollary~\ref{gh}.
\end{proof}

\begin{corollary}\label{trczs}
The toral rank conjecture holds for the restricted torus action on
$\zs$, that is, $\hrk\zs\ge 2^{\trk\zs}$.
\end{corollary}
\begin{proof}
We have that $\trk\zs=\trk\mathcal Z_{\mathcal K_\mathcal S}=m-n$
by Corollary~\ref{trkzs}, and $\hrk\mathcal Z_{\mathcal K_\mathcal
S}\ge2^{m-n}$ by~\cite[Cor.~1.4]{ca-lu09} or~\cite[\S3]{usti09}.
Therefore,
\[
  \hrk\zs\ge\hrk\mathcal Z_{\mathcal K_\mathcal S}\ge2^{m-n},
\]
as claimed.
\end{proof}

\begin{remark}
In fact, according to~\cite[Th.~3.2]{usti09}, the sharper bound
$\hrk\zs\ge2^{m-\mrk\sS}$ holds, where $\mrk\sS$ is the minimal
rank of maximal elements in~$\sS$. It equals $n$ (the rank of
$\sS$) if and only if $\sS$ is pure, that is, all maximal elements
of $\sS$ have the same rank.
\end{remark}

\end{document}